\documentclass[11pt]{article}
\usepackage{amssymb}
\usepackage{amsmath}
\usepackage{amsfonts}
\usepackage{amsthm}
\usepackage{enumerate}
\input amssym.def

\title{\bf Well-posedness and ill-posedness results for dissipative Benjamin-Ono equations}
\author{St\'ephane Vento, \\Universit\'e Paris-Est, \\Laboratoire d'Analyse
et de Math\'ematiques Appliqu\'ees,\\ 5 bd. Descartes, Cit\'e
Descartes, Champs-Sur-Marne,\\ 77454 Marne-La-Vall\'ee Cedex 2,
France}

\date{E-mail:\, stephane.vento@univ-paris-est.fr}

\numberwithin{equation}{section}

\newtheorem{theorem}{Theorem}[section]
\newtheorem{lemma}{Lemma}[section]
\newtheorem{proposition}{Proposition}[section]

\newtheorem{remark}{Remark}[section]

\def\R{\mathbb{R}}

\def\Z{\mathbb{Z}}
\def\C{\mathcal{C}}
\def\S{\mathcal{S}}
\def\H{\mathcal{H}}

\def\F{\mathcal{F}}
\def\eps{\varepsilon}
\def\supp{\mathop{\rm supp}\nolimits}
\def\sgn{\mathop{\rm sgn}\nolimits}
\def\pv{\mathop{\rm pv}\nolimits}
\newcommand\cro[1]{\langle #1 \rangle}

\begin{document}
\maketitle

\noindent {\bf Abstract.}\, We study the Cauchy problem for the
dissipative Benjamin-Ono equations $u_t+\H u_{xx}+|D|^\alpha
u+uu_x=0$ with $0\leq\alpha\leq 2$. When $0\leq\alpha< 1$, we show
the ill-posedness in $H^s(\R)$, $s\in\R$, in the sense that the
flow map $u_0\mapsto u$ (if it exists) fails to be $\C^2$ at the
origin. For $1<\alpha\leq 2$, we prove the global well-posedness
in $H^s(\R)$, $s>-\alpha/4$. It turns out that this index is
optimal.
\\


\section{Introduction, main results and notations}
\subsection{Introduction}
In this work we consider the Cauchy problem for the following
dissipative Benjamin-Ono equations
\begin{equation}\label{dBO}\tag{dBO}\left\{\begin{array}{ll}u_t+\H u_{xx}+|D|^\alpha u+uu_x=0,
\quad t>0, x\in\mathbb{R},\\u(0,\cdot)=u_0\in H^s(\R),\end{array}\right.\end{equation} with $0\leq\alpha\leq 2$, and where $\H$
is the Hilbert transform defined by
\[\mathcal{H}f(x)=\frac{1}{\pi}\pv\Big(\frac{1}{x}\ast
f\Big)(x)=\mathcal{F}^{-1}\big(-i\ \sgn(\xi)\hat{f}(\xi)\big)(x)\]
and $|D|^\alpha$ is the Fourier multiplier with symbol
$|\xi|^\alpha$.

\vskip 0.5cm

When $\alpha=0$, (\ref{dBO}) is the ordinary Benjamin-Ono equation
derived by Benjamin \cite{1967JFM....29..559B} and later by  Ono
\cite{MR0398275} as a model for one-dimensional waves in deep
water. The Cauchy problem for the Benjamin-Ono equation has been
extensively studied these last years. It has been proved in
\cite{MR533234} that (BO) is globally well-posed  in $H^s(\R)$ for
$s\geq 3$, and then for $s\geq 3/2$ in \cite{MR1097916} and
\cite{MR847994}. In \cite{MR2052470}, Tao
 get the well-posedness of this equation for
$s\geq 1$ by using a gauge transformation (which is a modified
version of the Cole-Hopf transformation). Recently, combining a
gauge transformation together with a Bourgain's method, Ionescu
and Kenig \cite{MR2291918} shown that one could go down to
$L^2(\R)$, and this seems to be, in some sense, optimal. It is
worth noticing that all these results have been obtained by
compactness methods. On the other hand,  Molinet, Saut and
 Tzvetkov \cite{MR1885293}  proved that, for all $s\in\R$, the
flow map $u_0\mapsto u$ is not of class $\mathcal{C}^2$ from
$H^s(\R)$ to $H^s(\R)$. Furthermore, building suitable families of
approximate solutions,   Koch and  Tzvetkov   proved in
\cite{MR2172940} that the flow map is actually not even uniformly
continuous on bounded sets of $H^s(\R)$, $s>0$. As an important
consequence of this, since a Picard iteration scheme would imply
smooth dependance upon the initial data, we see that such a scheme
cannot be used to get solutions in any space continuously embedded
in $\mathcal{C}([0,T];H^s(\R))$.

\vskip 0.5cm

When $\alpha=2$, (\ref{dBO}) is the so-called Benjamin-Ono-Burgers
equation
\begin{equation}\label{BOB}\tag{BOB}u_t+(\H-1)u_{xx}+uu_x=0.\end{equation} Edwin
and Robert \cite{MR830421} have derived (\ref{BOB}) by means of
formal asymptotic expansions in order to describe wave motions by
intense magnetic flux tube in the solar atmosphere. The
dissipative effects in that context are due to heat conduction.
(\ref{BOB}) has been studied in many papers, see
\cite{MR1101240,MR1811951,MR1722827}. Working in Bourgain's spaces
containing both dispersive and dissipative
effects\footnotemark[1], \footnotetext[1]{Such spaces were first
introduce by Molinet and Ribaud in \cite{MR1918236} for the
KdV-Burgers equation.} Otani showed in \cite{MR2174979} that
(\ref{BOB}) is globally well-posed in $H^s(\R)$, $s>-1/2$. In this
paper, we prove that this index is in fact critical  since the
flow map $u_0\mapsto u$ is not of class $\C^3$ from $H^s(\R)$ to
$H^s(\R)$, $s<-1/2$. Intriguingly, this index coincides with the
critical Sobolev space for the Burgers equation
$$u_t-u_{xx}+uu_x=0,$$ see \cite{MR1382829,MR1409926}. This result
is in a marked contrast with what occurs for the KdV-Burgers
equation which is well-posed above $H^{-1}(\R)$, see
\cite{MR1918236}.

\vskip 0.5 cm

Now consider the general case $0\leq\alpha\leq 2$. By running the
approach of \cite{MR1918236} combined with the smoothing relation
obtained in \cite{MR2174979}, we can only get that the problem
(\ref{dBO}) is well-posed in $H^s(\R)$ for $3/2<\alpha\leq 2$ and
$s>1/2-\alpha/2$. This was done by Otani in \cite{MR2303557}. Here
we improve this result by showing that (\ref{dBO}) is globally
well-posed in $H^s(\R)$, for $1<\alpha\leq 2$ and $s>-\alpha/4$.
It is worth comparing (\ref{dBO}) with the pure dissipative
equation \begin{equation}\label{dissip}u_t+|D|^\alpha
u+uu_x=0.\end{equation} In the Appendix, we show that
(\ref{dissip}) with $1<\alpha\leq 2$ is well-posed in $H^s(\R)$ as
soon as $s>3/2-\alpha$. The techniques we use are very common in
the context of semilinear parabolic problems and can be easily
adapted to (\ref{dBO}). In particular when $\alpha=2$, this
provides an alternative (and simpler) proof of our main result.
When $\alpha<2$, clearly we see that the dispersive part in
(\ref{dBO}) plays a key role in the low regularity of the
solution.

\vskip 0.5cm

We are going to perform a fixed point argument on the integral
formulation of (\ref{dBO}) in the weighted Sobolev space
\begin{equation}\label{Xbs}\|u\|_{X^{b,s}_\alpha}
=\|\cro{i(\tau-\xi|\xi|)+|\xi|^\alpha}^b\cro{\xi}^s\F
u(\tau,\xi)\|_{L^2(\R^2)}.\end{equation} This will be achieved by
deriving a bilinear estimate in these spaces. By Plancherel's
theorem and duality, it reduces to estimating a weighted
convolution of $L^2$ functions. In some regions where the
dispersive effect is too weak to recover the lost derivative in
the nonlinear term at low regularity ($s>-\alpha/4$), in
particular when considering the high-high interactions, we are led
to use a dyadic approach. In \cite{MR1854113}, Tao systematically
studied some nonlinear dispersive equations like KdV,
Schr\"{o}dinger or wave equation by using such dyadic
decomposition and orthogonality. Following the spirit of Tao's
works, we shall prove some estimates on dyadic blocks, which may
be of independent interest. Indeed, we believe that they could
certainly be used for other equations based on a Benjamin-Ono-type
dispersion.

\vskip 0.5cm

Next, we show that our well-posedness results turn out to be
sharp. Adapting the arguments used in \cite{MR1885293} to prove
the ill-posedness of (BO), we find that the solution map
$u_0\mapsto u$ (if it exists) cannot be $\C^3$ at the origin from
$H^s(\R)$ to $H^s(\R)$ as soon as $s<-\alpha/4$. See also
\cite{MR1466164,MR1918236,MR2038121,vento-2007} for situations
where this method applies. Note that we need to prove the
discontinuity of the third iterative term to obtain the condition
$s<-\alpha/4$, whereas  the second iterate is usually sufficient
to get an optimal result. On the other hand, we prove using
similar arguments, that in the case $0\leq\alpha< 1$, the solution
map fails to be $\C^2$ in any $H^s(\R)$, $s\in\R$. This is mainly
due to the fact that the operator $|D|^\alpha$ is too weak to
counterbalance the lost derivative which appears in the nonlinear
term $\partial_x u^2$.

\subsection{Main results}

Let us now formally state our results.

\begin{theorem}\label{main} Let $1<\alpha\leq 2$ and $u_0\in H^s(\R)$ with
$s>-\alpha/4$. Then for any $T>0$, there exists a unique solution
$u$ of (\ref{dBO}) in $$Z_T=\mathcal{C}([0,T]; H^s(\R))\cap
X^{1/2,s}_{\alpha,T}.$$ Moreover, the map $u_0\mapsto u$ is smooth
from $H^s(\R)$ to $Z_T$ and $u$ belongs to $\mathcal{C}((0,T],
H^\infty(\R))$.
\end{theorem}

\begin{remark} The spaces $X_{\alpha,T}^{b,s}$ are restricted
versions of $X_\alpha^{b,s}$ defined by the norm (\ref{Xbs}). See
Section \ref{sec-not} for a precise definition.
\end{remark}

\begin{remark} In \cite{MR2303557}, Otani studied a larger family
of dispersive-dissipative equations taking the form
\begin{equation}\label{gbob}u_t-|D|^{1+a}u_x+|D|^\alpha u+uu_x=0\end{equation} with $a\geq 0$ and $\alpha>0$.
He showed that (\ref{gbob}) is globally well-posed in $H^s(\R)$
provided $a+\alpha\leq 3$, $\alpha>(3-a)/2$ and
$s>-(a+\alpha-1)/2$. If $a=0$, it is clear that we get a better
result, at least when $\alpha< 2$. It will be an interesting
challenge to adapt our method of proofs to (\ref{gbob}) in the
case $a>0$.
\end{remark}

\begin{remark} Another interesting problem should be to consider
the periodic dissipative BO equations
\begin{equation}\label{dBOp}\left\{\begin{array}{ll}u_t+\H u_{xx}+|D|^\alpha u+uu_x=0,
\quad t>0, x\in\mathbb{T},\\u(0,\cdot)=u_0\in
H^s(\mathbb{T}),\end{array}\right.\end{equation} Recall that in
\cite{molinet-2006}, Molinet proved the global well-posedness of
the periodic BO equation in $L^2(\mathbb{T})$. To our knowledge,
equation (\ref{dBOp}) in the case $\alpha>0$ has never been
investigated.
\end{remark}

\vskip 0.5cm

Theorem \ref{main} is sharp in the following sense.

\begin{theorem}\label{th-12} Let $1\leq\alpha\leq 2$ and $s<-\alpha/4$. There does not exist $T>0$ such that
the Cauchy problem (\ref{dBO}) admits a unique local solution
defined on the interval $[0,T]$ and such that the flow map
$u_0\mapsto u$ is of class $\mathcal{C}^3$ in a neighborhood of
the origin from $H^s(\mathbb{R})$ to $H^s(\mathbb{R})$.
\end{theorem}

\vskip 0.5cm

In the case $0\leq\alpha< 1$, we have the following ill-posedness
result.

\begin{theorem}\label{th-01} Let $0\leq\alpha< 1$ and $s\in\R$. There does not exist $T>0$ such that
the Cauchy problem (\ref{dBO}) admits a unique local solution
defined on the interval $[0,T]$ and such that the flow map
$u_0\mapsto u$ is of class $\mathcal{C}^2$ in a neighborhood of
the origin from $H^s(\mathbb{R})$ to $H^s(\mathbb{R})$.
\end{theorem}

\begin{remark} At the end-point $\alpha=1$, our proof of Theorem \ref{th-01}
fails. However, Theorem \ref{th-12} provides the ill-posedness in
$H^s(\R)$, for $s<-1/4$. So, it is still not clear of what happens
to (\ref{dBO}) when $\alpha=1$ and $s\geq -1/4$.
\end{remark}

\vskip 0.5cm

The structure of our paper is as follows. We introduce a few
notation in the rest of this section. In Section \ref{sec-lin}, we
recall some estimates related to the linear (\ref{dBO}) equations.
Next, we prove the crucial bilinear estimate in Section
\ref{sec-bil}, which leads to the proof of Theorem \ref{main} in
Section \ref{sec-main}. Section \ref{sec-ill} is devoted to the
ill-posedness results (Theorems \ref{th-12} and \ref{th-01}).
Finally, we briefly study the dissipative equation (\ref{dissip})
in the Appendix.

\subsection{Notations}\label{sec-not}
When writing $A\lesssim B$ (for $A$ and $B$ nonnegative), we mean
that there exists $C>0$ independent of $A$ and $B$ such that
$A\leq CB$. Similarly define $A\gtrsim B$ and $A\sim B$. If
$A\subset\R^N$, $|A|$ denotes its Lebesgue measure and $\chi_A$
its characteristic function. For $f\in \S'(\R^N)$, we define its
Fourier transform $\F(f)$
 (or $\widehat{f}$) by
$$\F f(\xi)=\int_{\R^N}e^{-i\langle x,\xi\rangle}f(x)dx.$$
The Lebesgue spaces are endowed with the norm
$$\|f\|_{L^p(\R^N)}=\Big(\int_{\R^N}|f(x)|^pdx\Big)^{1/p},\quad
1\leq p<\infty$$ with the usual modification for $p=\infty$. We
also consider the space-time Lebesgue spaces $L^p_xL^q_t$ defined
by
$$\|f\|_{L^p_xL^q_t}=\Big\|\|f\|_{L^q_t(\R)}\Big\|_{L^p_x(\R)}.$$
For $b,s\in\R$, we define the Sobolev spaces $H^s(\R)$ and their
space-time versions $H^{b,s}(\R^2)$ by the norms
$$\|f\|_{H^s}=\Big(\int_\R\langle\xi\rangle^{2s}|\widehat{f}(\xi)|^2d\xi\Big)^{1/2},$$
$$\|u\|_{H^{b,s}}=\Big(\int_{\R^2}\langle\tau\rangle^{2b}\langle\xi\rangle^{2s}|\widehat{u}(\tau,\xi)|^2d\tau
d\xi\Big)^{1/2},$$ with $\langle\cdot\rangle=(1+|\cdot|^2)^{1/2}$.
Let $V(\cdot)$ be the free linear group associated to the linear
Benjamin-Ono equation, i.e.
$$\forall t\in\R,\
\F_x(V(t)\varphi)(\xi)=\exp(it\xi|\xi|)\widehat{\varphi}(\xi),\quad
\varphi\in\S'.$$ We will mainly work in the $X^{b,s}_\alpha$ space defined in (\ref{Xbs}), and in its restricted version
$X^{b,s}_{\alpha,T}$, $T\geq 0$, equipped with the norm
$$\|u\|_{X^{b,s}_{\alpha,T}}=\inf_{w\in
X^{b,s}_\alpha}\{\|w\|_{X^{b,s}_\alpha},\ w(t)=u(t)\textrm{ on }
[0,T]\}.$$ Note that since
$\F(V(-t)u)(\tau,\xi)=\widehat{u}(\tau+\xi|\xi|,\xi)$, we can
re-express the norm of $X^{b,s}_\alpha$ as \begin{eqnarray*}
\|u\|_{X^{b,s}_\alpha} &=& \big\|\langle
i\tau+|\xi|^{\alpha}\rangle^b\langle\xi\rangle^s\widehat{u}(\tau+\xi|\xi|,\xi)\big\|_{L^2(\R^2)}\\
&=& \big\|\langle
i\tau+|\xi|^{\alpha}\rangle^b\langle\xi\rangle^s\F(V(-t)u)(\tau,\xi)\big\|_{L^2(\R^2)}\\
&\sim & \|V(-t)u\|_{H^{b,s}}+\|u\|_{L^2_tH^{s+\alpha b}_x}.
\end{eqnarray*}
Finally, we denote by $S_\alpha$ the semigroup associated with the
free evolution of (\ref{dBO}), $$\forall t\geq 0,\
\F_x(S_\alpha(t)\varphi)(\xi)=\exp[it\xi|\xi|-|\xi|^{\alpha}t]\widehat{\varphi}(\xi),\
\varphi\in\S',$$ and we extend $S_\alpha$ to a linear operator
defined on the whole real axis by setting
\begin{equation}\label{salpha}\forall t\in\R,\
\F_x(S_\alpha(t)\varphi)(\xi)=\exp[it\xi|\xi|-|\xi|^{\alpha}|t|]\widehat{\varphi}(\xi),\
\varphi\in\S'.\end{equation}

\section{Linear estimates} \label{sec-lin}
In this section, we collect together several linear estimates on
the operators $S_\alpha$ introduced in (\ref{salpha}) and
$L_\alpha$ defined by $$L_\alpha: f\mapsto
\chi_{\R_+}(t)\psi(t)\int_0^tS_\alpha(t-t')f(t')dt'.$$

Recall that (\ref{dBO}) is equivalent to its integral formulation
\begin{equation}\label{eqint}u(t)=S_\alpha(t)u_0-\frac
12\int_0^tS_\alpha(t-t')\partial_x(u^2(t'))dt'.\end{equation} It
will be convenient to replace the local-in-time integral equation
(\ref{eqint}) with a global-in-time truncated integral equation.
Let $\psi$ be a cutoff function such that
$$\psi\in\mathcal{C}^\infty_0(\R),\quad \supp \psi\subset
[-2,2],\quad \psi\equiv 1\textrm{ on }[-1,1],$$ and define
$\psi_T(\cdot)=\psi(\cdot/T)$ for all $T>0$. We can replace
(\ref{eqint}) on the time interval $[0,T]$, $T<1$ by the equation
\begin{equation}\label{int}u(t)=\psi(t)\Big[S_\alpha(t)u_0-\frac{\chi_{\R_+}(t)}{2}\int_0^tS_\alpha(t-t')\partial_x(\psi_T^2(t')
u^2(t'))dt'\Big].
\end{equation}

Proofs of the results stated here can be obtained by a slight
modification of the linear estimates derived in \cite{MR1918236}.

\begin{lemma} For all $s\in\R$ and all $\varphi\in H^s(\R)$,
\begin{equation}\label{lin-free}\|\psi(t)S_\alpha(t)\varphi\|_{X^{1/2,s}_\alpha}\lesssim
\|\varphi\|_{H^s}.\end{equation}
\end{lemma}

\begin{lemma} Let $s\in\R$.
For all $0<\delta<1/2$ and all $v\in X^{-1/2+\delta,s}_\alpha$,
\begin{equation}\label{lin-for}\Big\|\chi_{\R_+}(t)\psi(t)\int_0^tS_\alpha(t-t')v(t')dt'\Big\|_{X^{1/2,s}_\alpha}\lesssim
\|v\|_{X^{-1/2+\delta,s}_\alpha}.\end{equation}
\end{lemma}

To globalize our solution, we will need the next lemma.
\begin{lemma}\label{glob} Let $s\in\R$ and $\delta>0$. Then for any
$f\in X^{-1/2+\delta,s}_\alpha$, $$t\longmapsto
\int_0^tS_\alpha(t-t')f(t')dt'\in\mathcal{C}(\R_+;H^{s+\alpha\delta}).$$
Moreover, if $(f_n)$ is a sequence satisfying $f_n\rightarrow 0$
in  $X^{-1/2+\delta,s}_\alpha$, then
$$\Big\|\int_0^tS_\alpha(t-t')f_n(t')dt'\Big\|_{L^\infty(\R_+;H^{s+\alpha\delta})}\longrightarrow
0.$$
\end{lemma}

\section{Bilinear estimates}\label{sec-bil}
\subsection{Dyadic blocks estimates}
\label{subsec-tao} We introduce Tao's $[k;Z]$-multipliers theory
\cite{MR1854113} and derive the dyadic blocks estimates for the
Benjamin-Ono equation.

Let $Z$ be any abelian additive group with an invariant measure
$d\eta$. For any integer $k\geq 2$ we define the hyperplane
$$\Gamma_k(Z)=\{(\eta_1,...,\eta_k)\in Z^k : \eta_1+...+\eta_k=0\}$$
which is endowed with the measure
$$\int_{\Gamma_k(Z)}f=\int_{Z^{k-1}}f(\eta_1,...,\eta_{k-1},-(\eta_1+...+\eta_{k-1}))d\eta_1...d\eta_{k-1}.$$
A $[k;Z]$-multiplier is defined to be any function
$m:\Gamma_k(Z)\rightarrow \mathbb{C}$. The multiplier norm
$\|m\|_{[k;Z]}$ is defined to be the best constant such that the
inequality
\begin{equation}\label{def-norm-mult}\Big|\int_{\Gamma_k(Z)}m(\eta)\prod_{j=1}^kf_j(\eta_j)\Big|\leq
\|m\|_{[k;Z]}\prod_{j=1}^k\|f_j\|_{L^2(Z)}\end{equation} holds for
all test functions $f_1,...,f_k$ on $Z$. In other words,
$$\|m\|_{[k;Z]}=\sup_{\substack{f_j\in\S(Z)\\\|f_j\|_{L^2(Z)}\leq
1}}\Big|\int_{\Gamma_k(Z)}m(\eta)\prod_{j=1}^kf_j(\eta_j)\Big|.$$
In his paper \cite{MR1854113}, Tao used the following notations.
Capitalized variables $N_j$, $L_j$ ($j=1,...,k$) are presumed to
be dyadic, i.e. range over numbers of the form $2^\ell$,
$\ell\in\Z$. In this paper, we only consider the case $k=3$, which
corresponds to the quadratic nonlinearity in the equation. It will
be convenient to define the quantities $N_{max}\geq N_{med}\geq
N_{min}$ to be the maximum, median and minimum of $N_1,N_2,N_3$
respectively. Similarly, define $L_{max}\geq L_{med}\geq L_{min}$
whenever $L_1,L_2,L_3>0$. The quantities $N_j$ will measure the
magnitude of frequencies of our waves, while $L_j$ measures how
closely our waves approximate a free solution.

Here we consider $[3;\R\times\R]$-multipliers and we parameterize
$\R\times\R$ by $\eta=(\tau,\xi)$ endowed with the Lebesgue
measure $d\tau d\xi$. Define
$$h_j(\xi_j)=\xi_j|\xi_j|,\quad \lambda_j=\tau_j-h_j(\xi_j),\quad j=1,2,3,$$
and the resonance function
$$h(\xi)=h_1(\xi_1)+h_2(\xi_2)+h_3(\xi_3).$$
By a dyadic decomposition of the variables $\xi_j$, $\lambda_j$,
$h(\xi)$, we will be led to estimate \begin{equation}\label{X}
\|X_{N_1,N_2,N_3,H,L_1,L_2,L_3}\|_{[3;\R\times\R]}
\end{equation}
where
\begin{equation}\label{defX}X_{N_1,N_2,N_3,H,L_1,L_2,L_3}=\chi_{|h(\xi)|\sim
H}\prod_{j=1}^3\chi_{|\xi_j|\sim N_j}\chi_{|\lambda_j|\sim L_j}.\end{equation}
From the identities
\begin{equation}\label{sigmaxj}\xi_1+\xi_2+\xi_3=0\end{equation}
and
$$\lambda_1+\lambda_2+\lambda_3+h(\xi)=0$$ on the support of the
multiplier, we see that (\ref{defX}) vanishes unless
\begin{equation}\label{nmax}N_{max}\sim N_{med}\end{equation} and
\begin{equation}\label{lmax}L_{max}\sim
\max(H,L_{med}).\end{equation}

\begin{lemma} On the support of $X_{N_1,N_2,N_3,H,L_1,L_2,L_3}$,
one has \begin{equation}\label{heq}H\sim
N_{max}N_{min}.\end{equation}
\end{lemma}

\begin{proof} Recall that
$$h(\xi)=\xi_1|\xi_1|+\xi_2|\xi_2|+\xi_3|\xi_3|.$$ By symmetry,
we can assume $|\xi_3|\sim N_{min}$. This forces by
(\ref{sigmaxj})  $\xi_1\xi_2<0$. Suppose for example $\xi_1>0$ and
$\xi_2<0$ (the other case being similar). Then if $\xi_3>0$,
$$h(\xi)=\xi_1^2-\xi_2^2+\xi_3^2=\xi_1^2-(\xi_1+\xi_3)^2+\xi_3^2=-2\xi_1\xi_3$$
and in this case $|h(\xi)|\sim N_{max}N_{min}$. Now if $\xi_3<0$,
then
$$h(\xi)=\xi_1^2-\xi_2^2-\xi_3^2=(\xi_2+\xi_3)^2-\xi_2^2-\xi_3^2=2\xi_2\xi_3$$
and it follows again that $|h(\xi)|\sim N_{max}N_{min}$.
\end{proof}

We are now ready to state the fundamental dyadic blocks estimates
for the Benjamin-Ono equation.

\begin{proposition} Let $N_1,N_2,N_3, H, L_1,L_2, L_3>0$
satisfying (\ref{nmax}), (\ref{lmax}), (\ref{heq}).
\begin{enumerate}
\item In the high modulation case $L_{max}\sim L_{med}\gg H$,
    we have
    \begin{equation}\label{esthigh}(\ref{X})\lesssim
    L_{min}^{1/2}N_{min}^{1/2}.\end{equation}
\item In the low modulation case $L_{max}\sim H$,
\begin{enumerate}
\item ((++) coherence) if $N_{max}\sim N_{min}$, then
    \begin{equation}\label{++}(\ref{X})\lesssim
    L_{min}^{1/2}L_{med}^{1/4},\end{equation}
\item ((+-) coherence) if $N_2\sim N_3\gg N_1$ and $H\sim
    L_1\gtrsim L_2, L_3$, we have for any  $\gamma>0$
    \begin{equation}\label{+-}(\ref{X})\lesssim
    L_{min}^{1/2}\min(N_{min}^{1/2},N_{max}^{1/
    2-1/2\gamma}N_{min}^{-1/2\gamma}L_{med}^{1/2\gamma}).\end{equation}
    Similarly for permutations of the indexes $\{1,2,3\}$.
    \item In all other cases, the multiplier (\ref{defX})
        vanishes.
\end{enumerate}
\end{enumerate}
\end{proposition}

\begin{proof} First we consider the high modulation case
$L_{max}\sim L_{med}\gg H$. Suppose for the moment that $L_1\geq
L_2\geq L_3$ and $N_1\geq N_2\geq N_3$. By using the comparison
principle (Lemma 3.1 in \cite{MR1854113}), we have
$$(\ref{X})\lesssim \|\chi_{|\xi_3|\sim N_3}\chi_{|\lambda_3|\sim
L_3}\|_{[3;\R\times\R]}.$$ By Lemma 3.14 and Lemma 3.6 in \cite{MR1854113},
$$(\ref{X})\lesssim \Big\|\|\chi_{|\lambda_3|\sim
L_3}\|_{[3;\R]}\chi_{|\xi_3|\sim N_3}\Big\|_{[3;\R]}\lesssim
L_3^{1/2}N_3^{1/2}.$$ It is clear from symmetry that (\ref{esthigh}) holds for any choice of $L_j$ and $N_j$, $j=1,2,3$.

Now we turn to the low modulation case $H\sim L_{max}$. Suppose
for the moment that $N_1\geq N_2\geq N_3$. The $\xi_3$ variable is
currently localized to the annulus $\{|\xi_3|\sim N_3\}$. By a
finite partition of unity we can restrict it further to a ball
$\{|\xi_3-\xi_3^0|\ll N_3\}$ for some $|\xi_3^0|\sim N_3$. Then by
box localisation (Lemma 3.13 in \cite{MR1854113}) we may localize
$\xi_1$, $\xi_2$ similarly to regions $\{|\xi_1-\xi_1^0|\ll N_3\}$
and $\{|\xi_2-\xi_2^0|\ll N_3\}$ where $|\xi_j^0|\sim N_j$. We may
assume that $|\xi_1^0+\xi_2^0+\xi_3^0|\ll N_3$ since we have
$\xi_1+\xi_2+\xi_3=0$. We summarize this symmetrically as
$$(\ref{X})\lesssim \Big\|\chi_{|h(\xi)|\sim
H}\prod_{j=1}^3\chi_{|\xi_j-\xi_j^0|\ll
N_{min}}\chi_{|\lambda_j|\sim L_j}\Big\|_{[3;\R\times\R]}$$ for
some $\xi_j^0$ satisfying $$|\xi_j^0|\sim N_j\textrm{ for
}j=1,2,3;\quad |\xi_1^0+\xi_2^0+\xi_3^0|\ll N_{min}.$$ Without
loss of generality, we assume $L_1\geq L_2\geq L_3$.  By Lemma
3.6, Lemma 3.1 and Corollary 3.10 in \cite{MR1854113}, we get
\begin{align*}(\ref{X}) &\lesssim \Big\|\chi_{|h(\xi)|\sim
H}\prod_{j=2}^3\chi_{|\xi_j-\xi_j^0|\ll
N_{min}}\chi_{|\lambda_j|\sim L_j}\Big\|_{[3;\R\times\R]}\\
&\lesssim |\{(\tau_2,\xi_2) : |\xi_2-\xi_2^0|\ll N_{min}, |\tau_2-h_2(\xi_2)|\sim L_2,\\
&\quad\quad  |\xi-\xi_2-\xi_3^0|\ll N_{min}, |\tau-\tau_2-h_3(\xi-\xi_2)|\sim L_3\}|^{1/2}
\end{align*}
for some $(\tau,\xi)\in \R\times\R$. For fixed $\xi_2$, the set of
possible $\tau_2$  ranges in an interval of length $O(L_3)$ and
vanishes unless
$$h_2(\xi_2)+h_3(\xi-\xi_2)=\tau+O(L_2).$$
On the other hand, inequality $|\xi-\xi_2-\xi_3^0|\ll N_{min}$
implies $|\xi+\xi_1^0|\ll N_{min}$, hence $$(\ref{X})\lesssim
L_3^{1/2}|\Omega_\xi|^{1/2}$$ for some $\xi$ such that
$|\xi+\xi_1^0|\ll N_{min}$ (in particular $|\xi|\sim N_1$) and
with
$$\Omega_\xi=\{\xi_2 : |\xi_2-\xi_2^0|\ll N_{min},
h_2(\xi_2)+h_3(\xi-\xi_2)=\tau+O(L_2)\}.$$ Let us write
$\Omega_\xi=\Omega_\xi^1\cup\Omega_\xi^2$ with
\begin{eqnarray*} \Omega_\xi^1 &=& \{\xi_2\in\Omega_\xi :
\xi_2(\xi-\xi_2)>0\}\\ \Omega_\xi^2 &=& \{\xi_2\in \Omega_\xi : \xi_2(\xi-\xi_2)<0\}.
\end{eqnarray*}
We need only to consider the three cases $N_1\sim N_2\sim N_3$,
$N_2\sim N_3\gg N_1$ and $N_1\sim N_2\gg N_3$ (the case $N_1\sim
N_3\gg N_2$ follows by symmetry).

\vskip 0.3cm \noindent Estimate of $|\Omega_\xi^1|$ : In
$\Omega_\xi^1$ we can assume $\xi_2>0$ and $\xi-\xi_2>0$ (the
other case being similar). Then we have
$$h_2(\xi_2)+h_3(\xi-\xi_2)=\xi_2^2+(\xi-\xi_2)^2=2\left(\xi_2-\frac
\xi 2\right)^2+\frac {\xi^2} 2$$ and thus
\begin{equation}\label{esto1}2\left(\xi_2-\frac \xi 2\right)^2+\frac {\xi^2}
2=\tau+O(L_2).\end{equation} If $N_1\sim N_2\sim N_3$, we see from
(\ref{esto1}) that $\xi_2$ variable is contained in the union of
two intervals of length $O(L_2^{1/2})$ at worst. Therefore $|\Omega_\xi^1|\lesssim L_2^{1/2}$ in this case.\\
If $N_1\sim N_2\gg N_3$, then
\begin{align*}\Big|\Big(\xi_2-\frac\xi
2\Big)+\frac{\xi_1^0}2\Big| &\leq \Big|\xi_2-\xi_2^0-\frac{\xi+\xi_1^0}2 -\xi_3^0\Big|+|\xi_1^0+\xi_2^0+\xi_3^0|\\
&\leq |\xi_2-\xi_2^0|+\frac 12|\xi+\xi_1^0|+|\xi_3^0|+|\xi_1^0+\xi_2^0+\xi_3^0|\\ &\lesssim N_{3}
\end{align*}
and we get $|\xi_2-\frac\xi 2|\sim N_1$. From (\ref{esto1}), we
see that we must have $N_1^2=O(L_2)$, which is in contradiction
with  $L_2\lesssim L_1\sim N_{max}N_{min}$. We deduce that the
multiplier vanishes in this region.\\
If $N_2\sim N_3\gg N_1$, then we have obviously $|\xi_2-\frac \xi
2|\sim N_2$ and, in the same way, the multiplier vanishes.

\vskip 0.3cm \noindent Estimate of $|\Omega_\xi^2|$ : We can
assume $\xi_2>0$ and $\xi-\xi_2<0$. It follows that
\begin{equation}\label{esto2}h_2(\xi_2)+h_3(\xi-\xi_2)=\xi_2^2-(\xi-\xi_2)^2=2\xi\left(\xi_2-\frac
\xi 2\right)=\tau+O(L_2).\end{equation} If $N_1\sim N_2\sim N_3$, we see
from (\ref{esto2}) that $\xi_2$ variable is contained in the union
of two intervals of length $O(N_1^{-1}L_2)$ at worst. But we have $L_2\lesssim L_1\sim N_1^2$ and thus $|\Omega_\xi^2|\lesssim L_2^{1/2}$
in this region.\\
If $N_1\sim N_2\gg N_3$, we have $|\xi_2-\frac \xi 2|\sim N_1$ as
previously and thus $N_1^2=O(L_2)$, the multiplier vanishes.\\
If $N_2\sim N_3\gg N_1$, then $|\xi_2-\frac \xi 2|\sim N_2$ and
for any $\gamma>0$, we have $|\xi_2-\frac\xi 2|\sim
N_2^{1-\gamma}|\xi_2-\frac\xi 2|^\gamma$. Therefore we see from
(\ref{esto2}) that $\xi_2$ variable is contained in the union of
two intervals of length
$O(N_2^{1-1/\gamma}N_1^{-1/\gamma}L_2^{1/\gamma})$ at worst, and
from $|\xi_2-\xi_2^0|\ll N_{min}$ we see that
$|\Omega_\xi^2|\lesssim N_{min}^{1/2}$, and (\ref{+-}) follows.
\end{proof}

\subsection{Bilinear estimate}
In this section we prove the following crucial bilinear estimate.

\begin{theorem}\label{thbil}  Let $1<\alpha\leq 2$ and $s>-\alpha/4$. For all $T>0$,
there exist $\delta,\nu>0$ such that for all $u,v\in
X^{1/2,s}_\alpha$ with compact support (in time) in $[-T,+T]$,
\begin{equation}\label{bil}\|\partial_x(uv)\|_{X^{-1/2+\delta,s}_\alpha}\lesssim
T^\nu\|u\|_{X^{1/2,s}_\alpha}\|v\|_{X^{1/2,s}_\alpha}.\end{equation}
\end{theorem}

To get the required contraction factor $T^\nu$ in our estimates,
the next lemma is very useful (see \cite{MR2303557}).

\begin{lemma}\label{contract} Let $f\in L^2(\R^2)$ with compact
support (in time) in $[-T, +T]$. For any $\theta>0$, there exists
$\nu=\nu(\theta)>0$ such that
$$\left\|\mathcal{F}^{-1}\Big(\frac{\widehat{f}(\tau,\xi)}{\langle\tau-\xi|\xi|\rangle^\theta}\Big)\right\|_{L^2_{xt}}\lesssim
T^\nu\|f\|_{L^2_{xt}}.$$
\end{lemma}

\begin{proof}[Proof of Theorem \ref{thbil}]
By duality, Plancherel and Lemma \ref{contract}, it suffices to
show that
$$\left\|\frac{\xi_3\cro{\xi_3}^s\cro{\xi_1}^{-s}\cro{\xi_2}^{-s}}{\cro{|\lambda_1|+|\xi_1|^\alpha}^{1/2}
\cro{|\lambda_2|+|\xi_2|^\alpha}^{1/2}
\cro{|\lambda_3|+|\xi_3|^\alpha}^{1/2-\delta}}\right\|_{[3;\R\times\R]}\lesssim
1.$$ By dyadic decomposition of the variables $\xi_j$,
$\lambda_j$, $h(\xi)$, we may assume $|\xi_j|\sim N_j$,
$|\lambda_j|\sim L_j$ and $|h(\xi)|\sim H$. By the translation
invariance of the $[k,Z]$-multiplier norm, we can always restrict
our estimate on  $L_j\gtrsim 1$ and $N_{max}\gtrsim 1$. The
comparison principle and orthogonality reduce our estimate to show
that
\begin{multline}\label{lowmod}\sum_{N_{max}\sim N_{med}\sim
N}\sum_{L_1,L_2,L_3\gtrsim
1}\frac{N_3\cro{N_3}^s\cro{N_1}^{-s}\cro{N_2}^{-s}}{(L_1+\cro{N_1}^\alpha)^{1/2}
(L_2+\cro{N_2}^\alpha)^{1/2}(L_3+\cro{N_3}^\alpha)^{1/2-\delta}}\\
\times \|X_{N_1,N_2,N_3,L_{max},L_1,L_2,L_3}\|_{[3;\R\times\R]}
\end{multline}
and
\begin{multline}\label{highmod}
\sum_{N_{max}\sim N_{med}\sim N}\sum_{L_{max}\sim
L_{med}}\sum_{H\ll
L_{max}}\frac{N_3\cro{N_3}^s\cro{N_1}^{-s}\cro{N_2}^{-s}}{(L_1+\cro{N_1}^\alpha)^{1/2}
(L_2+\cro{N_2}^\alpha)^{1/2}(L_3+\cro{N_3}^\alpha)^{1/2-\delta}}\\
\times \|X_{N_1,N_2,N_3,H,L_1,L_2,L_3}\|_{[3;\R\times\R]}
\end{multline}
are bounded, for all $N\gtrsim 1$.

We first show that $(\ref{highmod})\lesssim 1$.  For $s>-1/2$, one
has
$$N_3\cro{N_3}^s\cro{N_1}^{-s}\cro{N_2}^{-s}\lesssim
\cro{N_{min}}^{-s}N_{max}$$ and we get from (\ref{esthigh}),
\begin{align*}
(\ref{highmod}) &\lesssim \sum_{N_{max}\sim N}
\sum_{L_{max}\gg
NN_{min}}\frac{\cro{N_{min}}^{-s}NL_{min}^{1/2}N_{min}^{1/2}}{L_{min}^{1/2}(L_{max}+N^\alpha)^{1/2-\delta}
(L_{max}+\cro{N_{min}}^\alpha)^{1/2-\delta}L_{max}^\delta}\\
&\lesssim \sum_{N_{min}>0}\frac{N_{min}^{1/2}\cro{N_{min}}^{-s}N}{(NN_{min}+N^\alpha)^{1/2-\delta}(NN_{min}+\cro{N_{min}}^\alpha)^{1/2-\delta}}.
\end{align*}
When $N_{min}\lesssim 1$, we get
\begin{align*}
(\ref{highmod}) &
\lesssim \sum_{N_{min}\lesssim 1}\frac{N_{min}^{1/2}N}{N^{\alpha/2-\alpha\delta}(NN_{min})^{1/2-\delta}}\\ &\lesssim \sum_{N_{min}\lesssim 1}N_{min}^{\delta}
N^{(1-\alpha)/2+\delta(\alpha+1)}\\ &\lesssim 1
\end{align*}
for $\delta\ll 1$ and $\alpha>1$. When $N_{min}\gtrsim 1$, then
\begin{align*}
(\ref{highmod}) &\lesssim \sum_{N_{min}\gtrsim 1}\frac{N_{min}^{1/2-s}N}{(NN_{min})^{1/2-\delta-\eps}N^{\alpha\eps}(NN_{min})^{1/2-\delta}}\\
&\lesssim \sum_{N_{min}\gtrsim 1}N_{min}^{-1/2-s+2\delta+\eps}N^{2\delta-\eps(\alpha-1)}\\ &\lesssim 1
\end{align*}
for $\eps=2\delta/(\alpha-1)>0$, $\delta\ll 1$ and $s>-1/2$.

Now we show that $(\ref{lowmod})\lesssim 1$. We first deal with
the contribution where (\ref{++}) holds. In this case $N_{min}\sim
N_{max}$ and we get
\begin{align*}
(\ref{lowmod}) &\lesssim \sum_{L_{max}\sim N^2}\frac{N^{1-s}L_{min}^{1/2}L_{med}^{1/4}}{L_{min}^{1/2}
(L_{med}+N^\alpha)^{1/2}(L_{max}+N^\alpha)^{1/2-2\delta}L_{max}^\delta}\\ &\lesssim \frac{N^{1-s}}{N^{\alpha/4}N^{1-4\delta}}\\ &\lesssim
N^{-s-\alpha/4+4\delta}\lesssim 1
\end{align*}
for $s>-\alpha/4$ and $\delta\ll 1$.

Now we consider the contribution where (\ref{+-}) applies. By
symmetry it suffices to treat the two cases
\begin{eqnarray*}
& N_1\sim N_2\gg N_3, & H\sim L_3\gtrsim L_1,L_2,\\ & N_2\sim N_3\gg N_1, & H\sim L_1\gtrsim L_2,L_3.
\end{eqnarray*}
In the first case, estimate (\ref{+-}) applied with $\gamma=1$
yields
$$(\ref{X})\lesssim
L_{min}^{1/2}\min(N_3^{1/2},N_3^{-1/2}L_{med}^{1/2})\lesssim
L_{min}^{1/2}N_3^{1/4}N_3^{-1/4}L_{med}^{1/4}\sim
L_{min}^{1/2}L_{med}^{1/4}$$ and thus
\begin{align*}
(\ref{lowmod}) &\lesssim \sum_{N_3>0}\sum_{L_{max}\sim NN_3}\frac{N_3\cro{N_3}^sN^{-2s}L_{min}^{1/2}L_{med}^{1/4}}
{L_{min}^{1/2}(L_{med}+N^\alpha)^{1/2}(L_{max}+\cro{N_{min}}^\alpha)^{1/2-2\delta}L_{max}^\delta}\\
&\lesssim \sum_{N_3>0}\frac{N_3\cro{N_3}^sN^{-2s}}{N^{\alpha/4}(NN_3)^{1/2-2\delta}}\\ &\lesssim \sum_{N_3>0}N_3^{1/2+2\delta}\cro{N_3}^s
N^{-2s-\alpha/4-1/2+2\delta}.
\end{align*}
Since $-2s-\alpha/4-1/2+2\delta<0$, we may write
\begin{align*}(\ref{lowmod}) &\lesssim
\sum_{N_3>0}N_3^{1/2+2\delta}\cro{N_3}^{-s-\alpha/4-1/2+2\delta}\\
&\lesssim \sum_{N_3\lesssim 1}N_3^{1/2+2\delta}+\sum_{N_3\gtrsim 1}N_3^{-s-\alpha/4+4\delta}\\ &\lesssim 1
\end{align*}
for $\delta\ll 1$ and $s>-\alpha/4$.

Finally consider the case $N_2\sim N_3\gg N_1$, $H\sim L_1\gtrsim
L_2,L_3$. Let $0<\gamma\ll 1$. If we assume $N_{min}^{1/2}\lesssim
N_{max}^{1/2-1/2\gamma}N_{min}^{-1/2\gamma}L_{med}^{1/2\gamma}$,
i.e. $L_{med}\gtrsim N_{max}^{1-\gamma}N_{min}^{1+\gamma}$, then
we get from (\ref{+-}) that
\begin{align*}(\ref{lowmod}) &\lesssim \sum_{N_1>0}\sum_{L_{max}\sim
NN_1}\frac{\cro{N_1}^{-s}NL_{min}^{1/2}N_1^{1/2}}{L_{min}^{1/2}(L_{med}+N^\alpha)^{1/2-\delta}L_{max}^{1/2-\delta}L_{max}^\delta}\\
&\lesssim \sum_{N_1>0}\frac{N_1^{1/2}\cro{N_1}^{-s}N}{(N^{1-\gamma}N_1^{1+\gamma}+N^\alpha)^{1/2-\delta}(NN_1)^{1/2-\delta}}\\
&\lesssim \sum_{N_1>0}\frac{N_1^{\delta}\cro{N_1}^{-s}N^{1/2+\delta}}{(N^{1-\gamma}N_1^{1+\gamma}+N^\alpha)^{1/2-\delta}}.
\end{align*}
If $N_1\lesssim 1$, then
$$(\ref{lowmod}) \lesssim \sum_{N_1\lesssim
1}N_1^{\delta}N^{(1-\alpha)/2+\delta(1+\alpha)}\lesssim
1$$ for $\delta\ll 1$ and $\alpha>1$. If $N_1\gtrsim 1$, then
\begin{align*}
(\ref{lowmod}) &\lesssim
\sum_{N_1\gtrsim 1}\frac{N_1^{-s+\delta}N^{1/2+\delta}}{(N^{1-\gamma}N_1^{1+\gamma})^{1/2-\delta-\eps}N^{\alpha\eps}}\\
&\lesssim \sum_{N_1\gtrsim 1}N_1^{-s-1/2+(1+\gamma)(\delta+\eps)+\delta-\gamma/2} N^{\gamma(1/2-\delta)+2\delta-\eps(\alpha-1+\gamma)}\\ &\lesssim 1
\end{align*}
for $\delta,\gamma\ll 1$, $s>-1/2$  and
$\eps=[2\delta+\gamma(1/2-\delta)]/(\alpha-1+\gamma)>0$.\\
If we assume $N_{min}^{1/2}\gtrsim
N_{max}^{1/2-1/2\gamma}N_{min}^{-1/2\gamma}L_{med}^{1/2\gamma}$,
i.e. $L_{med}\lesssim N_{max}^{1-\gamma}N_{min}^{1+\gamma}$, we
get
\begin{align*}
(\ref{lowmod}) &\lesssim \sum_{N_1>0}\sum_{L_{max}\sim NN_1}
\frac{\cro{N_1}^{-s}NL_{min}^{1/2}N^{1/2-1/2\gamma}N_1^{-1/2\gamma}L_{med}^{1/2\gamma}}
{L_{min}^{1/2}(L_{med}+N^\alpha)^{1/2-\delta}L_{max}^{1/2-\delta} L_{max}^\delta}\\
&\lesssim \sum_{N_1>0}\sum_{L_{med}\lesssim N^{1-\gamma}N_1^{1+\gamma}}
\frac{N_1^{-1/2\gamma-1/2+\delta}\cro{N_1}^{-s}N^{1-1/2\gamma+\delta}L_{med}^{1/2\gamma}}
{(L_{med}+N^\alpha)^{1/2-\delta}}.
\end{align*}
When $N_1\lesssim 1$, we have \begin{align*}(\ref{lowmod}) &
\lesssim \sum_{N_1\lesssim
1}N_1^{-1/2\gamma-1/2+\delta}N^{1-1/2\gamma+\delta}N^{-\alpha/2+\alpha\delta}(N^{1-\gamma}N_1^{1+\gamma})^{1/2\gamma}\\
&\lesssim \sum_{N_1\lesssim 1}N_1^\delta
N^{(1-\alpha)/2+\delta(1+\alpha)}\lesssim 1
\end{align*}
for $\delta\ll 1$ and $\alpha>1$. When $N_1\gtrsim 1$, then
\begin{align*}
(\ref{lowmod}) &\lesssim \sum_{N_1\gtrsim 1}N_1^{-s-1/2-1/2\gamma+\delta}N^{1-1/2\gamma+\delta}(N^{1-\gamma}N_1^{1+\gamma})^{1/2\gamma-1/2+\delta+\eps}
N^{-\alpha\eps}\\ &\lesssim \sum_{N_1\gtrsim 1}
N_1^{-s-1/2+(1+\gamma)(\delta+\eps)+\delta-\gamma/2} N^{\gamma(1/2-\delta)+2\delta-\eps(\alpha-1+\gamma)}\\ &\lesssim 1
\end{align*}
as previously.
 This completes the
proof of Theorem \ref{thbil}.
\end{proof}

\section{Proof of Theorem \ref{main}}\label{sec-main}
In this section, we sketch the proof of Theorem \ref{main} (see
for instance \cite{MR1918236} for the details).

Actually, local existence of a solution is a consequence of the
following modified version of Theorem \ref{thbil}.

\begin{proposition} Given $s_c^+>-\alpha/4$, there exist
$\nu,\delta>0$ such that for any $s\geq s_c^+$ and any $u,v\in
X^{1/2,s}_\alpha$ with compact support in $[-T,+T]$,
\begin{equation}\label{modbil}\|\partial_x(uv)\|_{X^{-1/2+\delta,s}_\alpha}\lesssim
T^\nu(\|u\|_{X^{1/2,s_c^+}_\alpha}\|v\|_{X^{1/2,s}_\alpha}+\|u\|_{X^{1/2,s}_\alpha}\|v\|_{X^{1/2,s_c^+}_\alpha}).\end{equation}
\end{proposition}
Estimate (\ref{modbil}) is obtained thanks to (\ref{bil}) and the
triangle inequality $$\forall s\geq s_c^+,\
\langle\xi\rangle^s\leq
\langle\xi\rangle^{s_c^+}\langle\xi_1\rangle^{s-s_c^+}+\langle\xi\rangle^{s_c^+}\langle\xi-\xi_1\rangle^{s-s_c^+}.$$
Let $u_0\in H^s(\R)$ with $s>-\alpha/4$. Define $F(u)$ as
$$F(u)=\psi(t)\Big[S_\alpha(t)u_0-\frac{\chi_{\R_+}(t)}{2}\int_0^tS_\alpha(t-t')\partial_x(\psi_T^2(t')
u^2(t'))dt'\Big].$$ We shall prove that for $T\ll 1$, $F$ is
contraction in a ball of the Banach space $$Z=\{u\in
X^{1/2,s}_\alpha :\
\|u\|_Z=\|u\|_{X^{1/2,s_c^+}_\alpha}+\gamma\|u\|_{X^{1/2,s}_\alpha}<+\infty\},$$
where $\gamma$ is defined for all nontrivial $\varphi$ by
$$\gamma=\frac{\|\varphi\|_{H^{s_c^+}}}{\|\varphi\|_{H^s}}.$$
Combining (\ref{lin-free}), (\ref{lin-for}) as well as
(\ref{modbil}), it is easy to derive that $$\|F(u)\|_Z\leq
C(\|u_0\|_{H^{s_c^+}}+\gamma\|u_0\|_{H^s})+CT^\nu\|u\|_Z^2$$ and
$$\|F(u)-F(v)\|_Z\leq CT^\nu\|u-v\|_Z\|u+v\|_Z$$ for some
$C,\nu>0$. Thus, taking $T=T(\|u_0\|_{H^{s_c^+}})$ small enough,
we deduce that $F$ is contractive on the ball of radius
$4C\|u_0\|_{H^{s_c^+}}$ in $Z$. This proves the existence of a
solution $u$ to $u=F(u)$ in $X^{1/2,s}_{\alpha,T}$.

Following similar arguments of \cite{MR1918236}, it is not too
difficult to see that if $u_1,u_2\in X^{1/2,s}_{\alpha, T}$ are
solutions of (\ref{int}) and $0<\delta<T/2$, then there exists
$\nu>0$ such that
$$\|u_1-u_2\|_{X^{1/2,s}_{\alpha,\delta}}\lesssim
T^\nu\big(\|u_1\|_{X^{1/2,s}_{\alpha,T}}+\|u_2\|_{X^{1/2,s}_{\alpha,T}}\big)\|u_1-u_2\|_{X^{1/2,s}_{\alpha,\delta}},$$
which leads to $u_1\equiv u_2$ on $[0,\delta]$, and then on
$[0,T]$ by iteration. This proves the uniqueness of the solution.

It is straightforward to check that
$S_\alpha(\cdot)u_0\in\C(\R_+;H^s(\R))\cap\C(\R_+^\ast;H^\infty(\R))$.
Then it follows from Theorem \ref{thbil}, Lemma \ref{glob} and the
local existence of the solution that $$u\in\C([0,T];
H^s(\R))\cap\C((0,T];H^{s+\alpha\delta}(\R))$$ for some
$T=T(\|u_0\|_{H^{s_c^+}})$. By induction, we have
$u\in\C((0,T];H^\infty(\R))$. Taking the $L^2$-scalar product of
(\ref{dBO}) with $u$, we obtain that
$t\mapsto\|u(t)\|_{H^{s_c^+}}$ is nonincreasing on $(0,T]$. Since
the existence time of the solution depends only on the norm
$\|u_0\|_{H^{s_c^+}}$, this implies that the solution can be
extended globally in time.

\section{Ill-posedness results}\label{sec-ill}
This section is devoted to the proof of Theorems \ref{th-12} and
\ref{th-01}. We adopt the notation $p(\xi)=\xi|\xi|$.

Assume that $u$ is a solution of (\ref{dBO}) such that the
solution map $u_0\mapsto u$ is of class $\C^k$ ($k=2$ or $k=3$) at
the origin from $H^s(\R)$ to $H^s(\R)$. The relation
$$F(u,\varphi):=u(t,\varphi)-S_\alpha(t)\varphi+\frac 12\int_0^tS_\alpha(t-t')\partial_x(u^2(t',\varphi))dt' \equiv 0$$
combined with implicit function theorem gives
\begin{eqnarray*}
u_1(t,x) & := & \frac{\partial u}{\partial\varphi}(t,x,0)[h]=S_\alpha(t)h\\
u_2(t,x) & := & \frac{\partial^2u}{\partial\varphi^2}(t,x,0)[h,h]=\int_0^tS_\alpha(t-t')\partial_x(u_1(t'))^2dt'\\
u_3(t,x) & := & \frac{\partial^3u}{\partial\varphi^3}(t,x,0)[h,h,h]=\int_0^tS_\alpha(t-t')\partial_x(u_1(t')u_2(t'))dt'\\
& \textrm{etc} &
\end{eqnarray*}
Since the solution map is $\C^k$, we must have
\begin{equation}\label{estcont}\|u_k(t)\|_{H^s}\lesssim \|h\|_{H^s}^k,\quad\forall h\in
H^s(\R).\end{equation} In the sequel, we will show that
(\ref{estcont}) fails in the case $0\leq\alpha< 1$ and $k=2$, and
in the case $1\leq\alpha\leq 2$, $k=3$ and $s>-\alpha/4$.

\subsection{The case $0\leq\alpha< 1$}
It suffices to show the following lemma.
\begin{lemma} Let $0\leq\alpha< 1$ and $s\in\R$. There exists a sequence of functions
$\{h_N\}\subset H^s(\R)$ such that for all $T>0$,
$$\|h_N\|_{H^s}\lesssim 1,$$ and
$$\lim_{N\rightarrow\infty}\sup_{[0,T]}\Big\|\int_0^tS_\alpha(t-t')\partial_x(S_\alpha(t')h_N)^2dt'\Big\|_{H^s}=+\infty.$$
\end{lemma}

\begin{proof} We define $h_N$ by its Fourier transform\footnotemark[1]
\footnotetext[1]{As noticed in \cite{MR1885293}, $h_N$ is not a
real-valued function but the analysis works as well for $\Re e\
h_N$ instead of $h_N$.}
$$\widehat{h_N}(\xi)=\gamma^{-1/2}\chi_{I_1}(\xi)+\gamma^{-1/2}N^{-s}\chi_{I_2}(\xi)$$
with $I_1=[\gamma/2,\gamma]$, $I_2=[N,N+\gamma]$ and $N\gg 1$,
$\gamma\ll N$ to be chosen later. Then it is clear that
$\|h_N\|_{H^s}\sim 1$. Computing the Fourier transform of $u_2(t)$
leads to
\begin{align*}\mathcal{F}_x(u_2(t))(\xi) &=
c\xi\int_0^te^{i(t-t')p(\xi)}e^{-(t-t')|\xi|^\alpha}(e^{it'p(\xi)}e^{-t'|\xi|^\alpha}\widehat{h_N})^{\ast 2}(\xi)dt'\\
&= c\xi e^{itp(\xi)}e^{-t|\xi|^\alpha}\int_\R\widehat{h_N}(\xi_1)\widehat{h_N}(\xi-\xi_1)\\ &\quad \times
\int_0^te^{it'(p(\xi_1)+p(\xi-\xi_1)-p(\xi))}e^{-t'(|\xi_1|^\alpha+|\xi-\xi_1|^\alpha
-|\xi|^\alpha)}dt'd\xi_1\\ &=c\xi e^{itp(\xi)}e^{-t|\xi|^\alpha}\int_\R\widehat{h_N}(\xi_1)\widehat{h_N}(\xi-\xi_1)\\ &\quad \times\frac{e^{it(
p(\xi_1)+p(\xi-\xi_1)-p(\xi))}e^{-t(|\xi_1|^\alpha+|\xi-\xi_1|^\alpha
-|\xi|^\alpha)}-1}{i(p(\xi_1)+p(\xi-\xi_1)-p(\xi))-(|\xi_1|^\alpha+|\xi-\xi_1|^\alpha
-|\xi|^\alpha)}d\xi_1.
\end{align*}
Set
$$\chi(\xi,\xi_1)=i(p(\xi_1)+p(\xi-\xi_1)-p(\xi))-(|\xi_1|^\alpha+|\xi-\xi_1|^\alpha
-|\xi|^\alpha).$$ By support considerations, we have $\|u_2(t)\|_{H^s}\geq
\|v_2(t)\|_{H^s}$ with
\begin{equation}\label{v2}\F_x(v_2(t))(\xi)=cN^{-s}\gamma^{-1}\xi
e^{itp(\xi)}e^{-t|\xi|^\alpha}
\int_{K_\xi}\frac{e^{t\chi(\xi,\xi_1)}-1}{\chi(\xi,\xi_1)}d\xi_1\end{equation}
and $$K_\xi= \{\xi_1 : \xi_1\in I_1, \xi-\xi_1\in I_2\}\cup\{\xi_1
: \xi_1\in I_2, \xi-\xi_1\in I_1\}.$$
 We easily see that if $\xi_1\in K_\xi$, then
$\xi\in [N+\gamma/2, N+2\gamma]$ and
$$p(\xi_1)+p(\xi-\xi_1)-p(\xi)=2\xi_1(\xi_1-\xi)\sim \gamma N,$$
$$|\xi_1|^\alpha+|\xi-\xi_1|^\alpha-|\xi|^\alpha\lesssim  N^\alpha.$$
We deduce that for $\gamma=N^{\alpha-1}\ll N$, we have
$|\chi(\xi,\xi_1)|\sim N^\alpha$. Now define
$$t_N=(N+2\gamma)^{-\alpha-\eps}\sim N^{-\alpha-\eps}$$ so that $e^{-t_N|\xi|^\alpha}\gtrsim 1$.
By a Taylor expansion of the exponential function,
\begin{equation}\label{devexp}\frac{e^{t_N\chi(\xi,\xi_1)}-1}{\chi(\xi,\xi_1)}=t_N+R(t_N,\xi,\xi_1)\end{equation}
with $$|R(t_N,\xi,\xi_1)|\lesssim \sum_{k\geq
2}\frac{t_N^k|\chi(\xi,\xi_1)|^{k-1}}{k!}\lesssim
N^{-\alpha-2\eps}.$$ Therefore the main contribution of
(\ref{devexp}) in (\ref{v2}) is given by $t_N$, and since
$|K_\xi|\sim \gamma$, it follows that
\begin{align*}|\F_x(v_2(t_N))(\xi)| &\gtrsim
N^{-s+1}\gamma^{-1}e^{-(N+2\gamma)^{-\eps}}\gamma
N^{-\alpha-\eps}\chi_{[N+\gamma/2,N+2\gamma]}(\xi)\\ & \gtrsim
N^{-s+1-\alpha-\eps}\chi_{[N+\gamma/2,N+2\gamma]}(\xi).\end{align*}
We get the lower bound for the $H^s$-norm of $u_2(t_N)$
$$\|u_2(t_N)\|_{H^s}\gtrsim
N^{-s+1-\alpha-\eps}\Big(\int_{N+\gamma/2}^{N+2\gamma}(1+|\xi|^2)^{s}d\xi\Big)^{1/2}
\sim N^{1-\alpha-\eps}\gamma^{1/2}\sim N^{(1-\alpha)/2-\eps},$$ which leads to
$$\lim_{N\rightarrow\infty}\sup_{[0,T]}\|u_2(t)\|_{H^s}=+\infty$$
for $\eps\ll 1$ and $\alpha<1$, as we claim.
\end{proof}

\subsection{The case $1\leq\alpha\leq 2$}
Let $1\leq\alpha\leq 2$ and $s<-\alpha/4$. As previously, it
suffices to find a suitable sequence $\{h_N\}$ such that
$\|h_N\|_{H^s}\lesssim 1$ and
$$\lim_{N\rightarrow\infty}\sup_{[0,T]}\|u_3(t)\|_{H^s}=+\infty.$$
 For this purpose, we define the real-valued function $h_N$ by
\begin{equation}\label{hn}\widehat{h_N}(\xi)=N^{-s}\gamma^{-1/2}(\chi_{I_N}(\xi)+\chi_{I_N}(-\xi))\end{equation}
with $I_N=[N,N+2\gamma]$, $N\gg 1$ and $\gamma\ll N$ to be chosen
later. We have
$$
\F_x(u_3(t))(\xi)=
c\xi\int_0^te^{i(t-t')p(\xi)}e^{-(t-t')|\xi|^\alpha}\F_x(S_\alpha(t')h_N)\ast\F_x(u_2(t'))(\xi)dt'$$
and
\begin{multline*}
\F_x(S_\alpha(t')h_N)\ast\F_x(u_2(t'))(\xi) = c\int_{\R^2}
\widehat{h_N}(\xi_1)\widehat{h_N}(\xi_2-\xi_1)\widehat{h_N}(\xi-\xi_2)\xi_2\\
\times
e^{it'(p(\xi-\xi_2)+p(\xi_2))}e^{-t'(|\xi-\xi_2|^\alpha+|\xi_2|^\alpha)}\frac{e^{t\chi(\xi_2,\xi_1)}-1}{\chi(\xi_2,\xi_1)}d\xi_1d\xi_2.
\end{multline*}
Hence, we can write $u_3=v_3-w_3$ with
\begin{align*}
\F_x(v_3(t))(\xi) &=c\xi e^{itp(\xi)}e^{-t|\xi|^\alpha}\int_{\R^2}
\widehat{h_N}(\xi_1)\widehat{h_N}(\xi_2-\xi_1)\widehat{h_N}(\xi-\xi_2)
\frac{\xi_2}{\chi(\xi_2,\xi_1)}\\ & \quad\times \int_0^te^{it'(p(\xi_1)+p(\xi_2-\xi_1)+p(\xi-\xi_2)-p(\xi))}
e^{-t(|\xi_1|^\alpha+|\xi_2-\xi_1|^\alpha+|\xi-\xi_2|^\alpha-|\xi|^\alpha)}dt'd\xi_1d\xi_2\\ &=
c\xi e^{itp(\xi)}e^{-t|\xi|^\alpha}\int_{\R^2}
\widehat{h_N}(\xi_1)\widehat{h_N}(\xi_2-\xi_1)\widehat{h_N}(\xi-\xi_2)
\frac{\xi_2}{\chi(\xi_2,\xi_1)} \frac{e^{t\lambda(\xi,\xi_1,\xi_2)}-1}{\lambda(\xi,\xi_1,\xi_2)}d\xi_1d\xi_2
\end{align*}
and
\begin{align*}
\F_x(w_3(t))(\xi) &=c\xi e^{itp(\xi)}e^{-t|\xi|^\alpha}\int_{\R^2}
\widehat{h_N}(\xi_1)\widehat{h_N}(\xi_2-\xi_1)\widehat{h_N}(\xi-\xi_2)
\frac{\xi_2}{\chi(\xi_2,\xi_1)}\\ & \quad\times \int_0^te^{t'\chi(\xi,\xi_2)}dt'd\xi_1d\xi_2\\ &=
c\xi e^{itp(\xi)}e^{-t|\xi|^\alpha}\int_{\R^2}
\widehat{h_N}(\xi_1)\widehat{h_N}(\xi_2-\xi_1)\widehat{h_N}(\xi-\xi_2)
\frac{\xi_2}{\chi(\xi_2,\xi_1)}\frac{e^{t\chi(\xi,\xi_2)}-1}{\chi(\xi,\xi_2)}d\xi_1d\xi_2
\end{align*}
where we set $$\lambda(\xi,\xi_1,\xi_2)=
i(p(\xi_1)+p(\xi_2-\xi_1)+p(\xi-\xi_2)-p(\xi))-(|\xi_1|^\alpha+|\xi_2-\xi_1|^\alpha+|\xi-\xi_2|^\alpha-|\xi|^\alpha).$$
Let $t_N=(N+4\gamma)^{-\alpha-\eps}$ for some $0<\eps\ll 1$. We
get
$$|\F_x(v_3(t_N))(\xi)|\chi_{[N+3\gamma,N+4\gamma]}(\xi) \gtrsim
N^{-3s+1}\gamma^{-3/2}\Big|\int_{K_\xi}\frac{\xi_2}{\chi(\xi_2,\xi_1)}\frac{e^{t_N\lambda(\xi,\xi_1,\xi_2)}-1}{\lambda(\xi,\xi_1,\xi_2)}d\xi_1,d\xi_2\Big|$$
where $K_\xi=K_\xi^1\cup K_\xi^2\cup K_\xi^3$ and
\begin{eqnarray*}
K_\xi^1 &=& \{(\xi_1,\xi_2) : \xi_1\in I_N, \xi_2-\xi_1\in I_N, \xi-\xi_2\in -I_N\},\\
K_\xi^2 &=& \{(\xi_1,\xi_2) : \xi_1\in I_N, \xi_2-\xi_1\in -I_N, \xi-\xi_2\in I_N\},\\
K_\xi^3 &=& \{(\xi_1,\xi_2) : \xi_1\in -I_N, \xi_2-\xi_1\in I_N, \xi-\xi_2\in I_N\}.
\end{eqnarray*}
If $\xi\in [N+3\gamma, N+4\gamma]$ and $(\xi_1,\xi_2)\in K_\xi$,
we easily see that
$$\Big|\frac{\xi_2}{\chi(\xi_2,\xi_1)}\Big|\sim N^{-1}$$ and
$$p(\xi_1)+p(\xi_2-\xi_1)+p(\xi-\xi_2)-p(\xi) \sim \gamma^2,$$
$$|\xi_1|^\alpha+|\xi_2-\xi_1|^\alpha+|\xi-\xi_2|^\alpha-|\xi|^\alpha\sim
N^\alpha.$$ Thus we are led to choose $\gamma= N^{\alpha/2}\ll N$
for $N\gg 1$ so that $|\lambda(\xi,\xi_1,\xi_2)|\sim N^\alpha$.
Then it follows that
$$\Big|\frac{e^{t_N\lambda(\xi,\xi_1,\xi_2)}-1}{\lambda(\xi,\xi_1,\xi_2)}\Big|=|t_N|+O(N^{-\alpha-2\eps}).$$
Consequently,
\begin{align*}|\F_x(v_3(t_N))(\xi)|\chi_{[N+3\gamma,N+4\gamma]}(\xi)
 &\gtrsim
N^{-3s+1}\gamma^{-3/2}N^{-1}\gamma^2N^{-\alpha-\eps}\chi_{[N+3\gamma,N+4\gamma]}(\xi)\\
&\sim
N^{-3s-\alpha-\eps}\gamma^{1/2}\chi_{[N+3\gamma,N+4\gamma]}(\xi)\\
&\sim N^{-3s-3\alpha/4-\eps}\chi_{[N+3\gamma,
N+4\gamma]}(\xi),\end{align*} since $|K_\xi|\sim \gamma^2$.

Concerning $w_3$, we verify that for $(\xi_1,\xi_2)\in K_\xi$, we
have $|\chi(\xi,\xi_2)|\gtrsim \gamma N$ and then
\begin{align*}|\F_x(w_3(t_N))(\xi)|\chi_{[N+3\gamma,
N+4\gamma]}(\xi) &
\lesssim N^{-3s+1}\gamma^{-3/2}\gamma^2N^{-1}(\gamma
N)^{-1}\chi_{[N+3\gamma, N+4\gamma]}(\xi)\\ &\sim N^{-3s-1}\gamma^{-1/2}\chi_{[N+3\gamma, N+4\gamma]}(\xi)\\ &\sim N^{-3s-1-\alpha/4}\chi_{[N+3\gamma, N+4\gamma]}(\xi)
\end{align*}
Since $-3s-1-\alpha/4<-3s-3\alpha/4-\eps$ for $\alpha<2$, we
deduce that the main contribution in the $H^s$-norm of $u_3$ is
given by $\|v_3\|_{H^s}$, that is,
$$\|u_3(t_N)\|_{H^s}\gtrsim N^{-3s-3\alpha/4-\eps}\gamma^{1/2}N^s\sim
N^{-2s-\alpha/2-\eps},$$ and we find the condition
$$-2s-\alpha/2>0,\quad\textrm{i.e.}\quad s<-\alpha/4.$$ When
$\alpha=2$, the contributions of $v_3$ and $w_3$ are equivalent,
and we must proceed with a bit more care, by considering directly
the difference $u_3=v_3-w_3$. More precisely, for $\gamma=\eps
N\ll N$, we have
$$|\lambda(\xi,\xi_1,\xi_2)|\sim |\chi(\xi,\xi_2)|\sim N^2.$$
Noticing that
$$\lambda(\xi,\xi_1,\xi_2)-\chi(\xi,\xi_2)=\chi(\xi_2,\xi_1),$$ we
deduce
$$\Big|\frac{e^{t_N\lambda(\xi,\xi_1,\xi_2)}-1}{\lambda(\xi,\xi_1,\xi_2)}-\frac{e^{t_N\chi(\xi,\xi_2)}-1}{\chi(\xi,\xi_2)}\Big|=t_N^2|\chi(\xi_2,\xi_1)|
+O(t_N^3N^2|\chi(\xi_2,\xi_1)|)
$$
Setting again $t_N=N^{-2-\eps}$, and since $|\xi_2|\sim N$, it
follows that
$$|\F_x(u_3(t_N))(\xi)|\chi_{[N+3\gamma,N+4\gamma]}\gtrsim
N^{-3s+1}\gamma^{-3/2}\gamma^2NN^{-4-2\eps}\chi_{[N+3\gamma,N+4\gamma]}(\xi)$$ and thus
$$\|u_3(t_N)\|_{H^s}\gtrsim N^{-2s-2-2\eps}\gamma\sim N^{-2s-1-2\eps},$$
which tends to infinity as soon as $-2s-1>0$, i.e. $s<-1/2$.

\appendix
\section{Appendix}
We prove here that the pure dissipative equation
\begin{equation}\label{gbur}u_t+|D|^\alpha u+uu_x=0\end{equation} for $1<\alpha\leq 2$ is well-posed in
$H^s(\R)$, $s>s_\alpha$ where
$$s_\alpha=\frac 32-\alpha,$$ and that the solution map fails to be
smooth when $s<s_\alpha$. The method of proof is classical and is
based on the smoothing properties of the generalized heat kernel
$$G_\alpha(t,x)=\frac{1}{2\pi}\int_\R e^{ix\xi}e^{-t|\xi|^\alpha}d\xi,\quad t>0.$$
\begin{theorem}\label{thgbur} Let $1<\alpha\leq 2$,
$s>s_\alpha$ and $u_0\in H^s(\R)$. Then there exist $T>0$ and a
unique solution $u\in\C([0,T]; H^s(\R))$ of (\ref{gbur}) such that
\begin{equation}\label{class1}\sup_{t\in[0,T]}\|u(t)\|_{H^s}<\infty\ \textrm{ if
}\
1<\alpha\leq 3/2,\end{equation}
\begin{equation}\label{class2}\sup_{t\in[0,T]}\|u(t)\|_{H^s}+\sup_{t\in[0,T]}t^\beta\|u(t)\|_{L^{2/(\alpha-1)}}<\infty\
\textrm{ if }\ 3/2<\alpha\leq 2\end{equation} where
$\beta=-s/\alpha+(2-\alpha)/2\alpha$. The flow map $u_0\mapsto u$
from $H^s(\R)$ into the class defined by
(\ref{class1})-(\ref{class2}) is locally Lipschitz. Moreover, if
$\|u_0\|_{H^s}$ is small enough, the solution can be extended to
any time interval.
\end{theorem}

\begin{proof} Observe that for any $p\in[1,\infty]$ and $\rho\geq
0$, we have \begin{equation}\label{estGt}\||D|^\rho
G_\alpha(t)\|_{L^p}=ct^{-(1-1/p)/\alpha-\rho/\alpha}.\end{equation}
We use the Picard iteration theorem to show that the map $F$
defined as $$F(u)=G_\alpha(t)\ast u_0-\frac
12\int_0^tG_\alpha(t-t')\ast\partial_xu^2(t')dt'$$ has a fixed
point in suitable Banach space.

We first consider the case $1<\alpha\leq 3/2$, and we choose
$s_\alpha<s<1/2$. Set $X_T=\C([0,T]; H^s(\R))$ endowed with the
norm $\|u\|_{X_T}=\sup_{[0,T]}\|u(t)\|_{H^s}$. By Young inequality
and (\ref{estGt}), we have
\begin{equation}\label{lin1}\|G_\alpha(t)\ast u_0\|_{H^s}\lesssim
\|G_\alpha(t)\|_{L^1}\|u_0\|_{H^s}\lesssim \|u_0\|_{H^s}.\end{equation} Using
the fractional Leibniz rule, we get
\begin{align*}
\int_0^t\|G_\alpha(t-t')\ast\partial_xu^2(t')\|_{H^s}dt' &\lesssim \int_0^t\|\partial_x G_\alpha(t-t')\|_{L^{(s+\frac 12)^{-1}}}\|\cro{D}^su^2(t')\|_{L^{1/(1-s)}}dt'\\
&\lesssim \int_0^t(t-t')^{s/\alpha-3/2\alpha}\|u(t')\|_{L^{(\frac 12-s)^{-1}}}\|u(t')\|_{H^s}dt'.
\end{align*}
Since $0<s<1/2$, we can take advantage of the Sobolev embedding
$H^s(\R)\hookrightarrow L^{(\frac 12-s)^{-1}}(\R).$ Since
$s/\alpha-3/2\alpha>-1$, we conclude
\begin{equation}\label{estint1}\int_0^t\|G_\alpha(t-t')\ast\partial_xu^2(t')\|_{H^s}dt'\lesssim
T^\nu\|u\|_{X_T}^2\end{equation} with $\nu=1+s/\alpha-3/2\alpha>0$. Gathering (\ref{lin1}) and (\ref{estint1}) we infer
$$\|F(u)\|_{X_T}\lesssim \|u_0\|_{H^s}+T^\nu\|u\|_{X_T}^2$$ and in
the same way, $$\|F(u)-F(v)\|_{X_T}\lesssim
T^\nu(\|u\|_{X_T}+\|v\|_{X_T})\|u-v\|_{X_T}.$$ This proves that
for $T\ll 1$, $F$ is contractive in a ball of $X_T$.

Now we solve (\ref{gbur}) in the case $3/2<\alpha\leq 2$ and
$s_\alpha<s<0$. Define
$Y_T=\C([0,T];H^s(\R))\cap\C^\beta([0,T];L^{2/(\alpha-1)}(\R))$
equipped with the norm
$$\|u\|_{Y_T}=\sup_{t\in[0,T]}\|u(t)\|_{H^s}+\sup_{t\in[0,T]}t^\beta\|u(t)\|_{L^{2/(\alpha-1)}}.$$
By Young inequality, we get $$\|G_\alpha(t)\ast
u_0\|_{L^{2/(\alpha-1)}}=\|\cro{D}^{-s}G_\alpha(t)\ast\cro{D}^su_0\|_{L^{2/(\alpha-1)}}\lesssim
\|\cro{D}^{-s}G_\alpha(t)\|_{L^{2/\alpha}}\|u_0\|_{H^s},$$ and it
follows from (\ref{estGt}) that
$$t^\beta\|\cro{D}^{-s}G_\alpha(t)\|_{L^{2/\alpha}}\lesssim
t^\beta(t^{-(2-\alpha)/2\alpha}+t^{-(2-\alpha)/2\alpha+s/\alpha})\lesssim
\cro{T}^{-s/\alpha}.$$ Now we deal with the nonlinear term. Using
the Sobolev embedding $L^{(\frac 12-s)^{-1}}(\R)\hookrightarrow
H^s(\R)$ valid for any $-1/2<s<0$, we obtain
\begin{align*}
\int_0^t\|G_\alpha(t-t')\ast\partial_x u^2(t')\|_{H^s}dt' &\lesssim \int_0^t\|\partial_xG_\alpha(t-t')\|_{L^{(\frac 52-s-\alpha)^{-1}}}
\|u^2(t')\|_{L^{1/(\alpha-1)}}dt'\\ &\lesssim \int_0^t(t-t')^{-s/\alpha-1+1/2\alpha}t'^{-2\beta}t'^{2\beta}\|u(t')\|_{L^{2/(\alpha-1)}}^2dt'\\ &\lesssim T^\nu
\|u\|_{Y_T}^2
\end{align*}
with $\nu=-s/\alpha+1/2\alpha-2\beta>0$. By similar calculations,
we get
\begin{align*}
t^\beta\int_0^t\|G_\alpha(t-t')\ast\partial_x u^2(t')\|_{L^{2/(\alpha-1)}}dt' &\lesssim t^\beta\int_0^t \|\partial_xG_\alpha(t-t')\|_{L^{2/(3-\alpha)}}
\|u^2(t')\|_{L^{1/(\alpha-1)}}dt'\\ &\lesssim t^\beta\int_0^t(t-t')^{-(\alpha+1)/2\alpha}t'^{-2\beta}dt' \|u\|_{Y_T}^2\\ &\lesssim T^\nu\|u\|_{Y_T}^2
\end{align*}
with $\nu=1-(\alpha+1)/2\alpha-\beta>0$. Finally, one has
$$\|F(u)\|_{Y_T}\lesssim
\cro{T}^\nu\|u_0\|_{H^s}+T^\nu\|u\|_{Y_T}^2$$ and the claim follows.
\end{proof}

\begin{remark}  Let
$U_\alpha(t)=\F_\xi^{-1}(e^{it\xi|\xi|}e^{-t|\xi|^\alpha})$ be the
fundamental solution of the linear (\ref{dBO}) equation. Using
that $|\F_xU_\alpha(t)|=|\F_xG_\alpha(t)|$ as well as the
well-known estimate $\|f\|_{L^p}\lesssim \|\hat{f}\|_{L^{p'}}$,
$p\geq 2$, $1/p+1/p'=1$, we easily check that Theorem \ref{thgbur}
holds for (\ref{dBO}) equation.
\end{remark}

Finally, we show that Theorem \ref{thgbur} is sharp.

\begin{theorem} Let $1<\alpha\leq 2$ and $s<s_\alpha$. The the
solution map $u_0\mapsto u$ associated with (\ref{gbur}) (if it
exists) is not of class $\C^2$ from $H^s(\R)$ to
$\C([0,T];H^s(\R))$.
\end{theorem}

\begin{proof} The proof is similar to that of Theorems
\ref{th-12} and \ref{th-01}. Define $h_N$ as in (\ref{hn}) and
consider the high-high interactions in the convolution product
$(e^{-t|\xi|^\alpha}h_N)\ast (e^{-t|\xi|^\alpha}h_N)$. We get that
for $\xi\in[2N, 2N+4\gamma]$, $\gamma=N^{1-\eps}$ and $t_N\sim
N^{-\alpha-\eps}$,
$$|\F_x(u_2(t_N))(\xi)|\gtrsim
N^{-2s-\alpha+1-\eps}\chi_{[2N,2N+4\gamma]}(\xi)$$ where $u_2$ is
defined by
$$u_2(t)=\int_0^tG_\alpha(t-t')\ast\partial_x(G_\alpha(t')\ast
h_N)^2dt'.$$ We conclude that
$$\|u_2(t_N)\|_{H^s}\gtrsim
N^{-s-\alpha+1-\eps}\gamma^{1/2}\gtrsim
N^{-s+3/2-\alpha-3\eps/2}\rightarrow +\infty$$ as soon as $s<3/2-\alpha$.
\end{proof}

\section*{Acknowledgments}
The author wishes to express his gratitude to Francis Ribaud for
his encouragement and precious advice.

\bibliographystyle{plain}
\bibliography{ref}

\begin{thebibliography}{10}

\bibitem{MR1409926}
D.~Bekiranov.
\newblock The initial-value problem for the generalized {B}urgers' equation.
\newblock {\em Differential Integral Equations}, 9(6):1253--1265, 1996.

\bibitem{1967JFM....29..559B}
T.~B. {Benjamin}.
\newblock {Internal waves of permanent form in fluids of great depth}.
\newblock {\em Journal of Fluid Mechanics}, 29:559--592, 1967.

\bibitem{MR1466164}
J.~Bourgain.
\newblock Periodic {K}orteweg de {V}ries equation with measures as initial
  data.
\newblock {\em Selecta Math. (N.S.)}, 3(2):115--159, 1997.

\bibitem{MR1101240}
D.~B. Dix.
\newblock Temporal asymptotic behavior of solutions of the
  {B}enjamin-{O}no-{B}urgers equation.
\newblock {\em J. Differential Equations}, 90(2):238--287, 1991.

\bibitem{MR1382829}
D.~B. Dix.
\newblock Nonuniqueness and uniqueness in the initial-value problem for
  {B}urgers' equation.
\newblock {\em SIAM J. Math. Anal.}, 27(3):708--724, 1996.

\bibitem{MR830421}
P.~M. Edwin and B.~Roberts.
\newblock The {B}enjamin-{O}no-{B}urgers equation: an application in solar
  physics.
\newblock {\em Wave Motion}, 8(2):151--158, 1986.

\bibitem{MR1811951}
A.~S. Fokas and L.~Luo.
\newblock Global solutions and their asymptotic behavior for
  {B}enjamin-{O}no-{B}urgers type equations.
\newblock {\em Differential Integral Equations}, 13(1-3):115--124, 2000.

\bibitem{MR2291918}
A.~D. Ionescu and C.~E. Kenig.
\newblock Global well-posedness of the {B}enjamin-{O}no equation in
  low-regularity spaces.
\newblock {\em J. Amer. Math. Soc.}, 20(3):753--798 (electronic), 2007.

\bibitem{MR847994}
R.~J. I{\'o}rio, Jr.
\newblock On the {C}auchy problem for the {B}enjamin-{O}no equation.
\newblock {\em Comm. Partial Differential Equations}, 11(10):1031--1081, 1986.

\bibitem{MR2172940}
H.~Koch and N.~Tzvetkov.
\newblock Nonlinear wave interactions for the {B}enjamin-{O}no equation.
\newblock {\em Int. Math. Res. Not.}, (30):1833--1847, 2005.

\bibitem{molinet-2006}
L.~Molinet.
\newblock Global well-posedness in {$L^2$} for the periodic {B}enjamin-{O}no
  equation, 2006.

\bibitem{MR1918236}
L.~Molinet and F.~Ribaud.
\newblock On the low regularity of the {K}orteweg-de {V}ries-{B}urgers
  equation.
\newblock {\em Int. Math. Res. Not.}, (37):1979--2005, 2002.

\bibitem{MR2038121}
L.~Molinet and F.~Ribaud.
\newblock Well-posedness results for the generalized {B}enjamin-{O}no equation
  with small initial data.
\newblock {\em J. Math. Pures Appl. (9)}, 83(2):277--311, 2004.

\bibitem{MR1885293}
L.~Molinet, J.-C. Saut, and N.~Tzvetkov.
\newblock Ill-posedness issues for the {B}enjamin-{O}no and related equations.
\newblock {\em SIAM J. Math. Anal.}, 33(4):982--988 (electronic), 2001.

\bibitem{MR0398275}
H.~Ono.
\newblock Algebraic solitary waves in stratified fluids.
\newblock {\em J. Phys. Soc. Japan}, 39(4):1082--1091, 1975.

\bibitem{MR2174979}
M.~Otani.
\newblock Bilinear estimates with applications to the generalized
  {B}enjamin-{O}no-{B}urgers equations.
\newblock {\em Differential Integral Equations}, 18(12):1397--1426, 2005.

\bibitem{MR2303557}
M.~Otani.
\newblock Well-posedness of the generalized {B}enjamin-{O}no-{B}urgers
  equations in {S}obolev spaces of negative order.
\newblock {\em Osaka J. Math.}, 43(4):935--965, 2006.

\bibitem{MR1097916}
G.~Ponce.
\newblock On the global well-posedness of the {B}enjamin-{O}no equation.
\newblock {\em Differential Integral Equations}, 4(3):527--542, 1991.

\bibitem{MR533234}
J.-C. Saut.
\newblock Sur quelques g\'en\'eralisations de l'\'equation de {K}orteweg-de
  {V}ries.
\newblock {\em J. Math. Pures Appl. (9)}, 58(1):21--61, 1979.

\bibitem{MR1854113}
T.~Tao.
\newblock Multilinear weighted convolution of {$L\sp 2$}-functions, and
  applications to nonlinear dispersive equations.
\newblock {\em Amer. J. Math.}, 123(5):839--908, 2001.

\bibitem{MR2052470}
T.~Tao.
\newblock Global well-posedness of the {B}enjamin-{O}no equation in {$H\sp
  1({\bf R})$}.
\newblock {\em J. Hyperbolic Differ. Equ.}, 1(1):27--49, 2004.

\bibitem{vento-2007}
S.~Vento.
\newblock Sharp well-posedness results for the generalized {B}enjamin-{O}no
  equation with high nonlinearity, 2007.

\bibitem{MR1722827}
L.~Zhang.
\newblock Local {L}ipschitz continuity of a nonlinear bounded operator induced
  by a generalized {B}enjamin-{O}no-{B}urgers equation.
\newblock {\em Nonlinear Anal.}, 39(3, Ser. A: Theory Methods):379--402, 2000.

\end{thebibliography}

\end{document}